\documentclass[twocolumn]{article}
\usepackage{amsfonts}
\usepackage{amsmath}
\usepackage{amssymb}
\usepackage{amsthm}
\usepackage{graphicx}

\usepackage{epsfig}
\usepackage{fancyhdr}
\usepackage{stmaryrd}
\usepackage{epsf}
\usepackage{amsmath,amsfonts,amssymb, amsbsy}
\usepackage{pifont}
\usepackage{latexsym}
\usepackage{epsfig}
\usepackage{rotating}
\usepackage{graphicx}
\usepackage{float}
\usepackage{pstricks, pst-node, pst-text, pst-all, multido}
\usepackage{times}
\usepackage{amsthm}
\usepackage{amsxtra}

\theoremstyle{definition}

\theoremstyle{remark}

\numberwithin{equation}{section}

\begin{document}

\title{NUMERICAL AND EXPERIMENTAL STUDY FOR A BEAM SYSTEM WITH LOCAL UNILATERAL CONTACT MODELING SATELLITE SOLAR ARRAYS}
\date{}

\author{Hamad Hazim$^{(1)}$, Bernard Rousselet$^{(1)}$ and Neil Ferguson$^{(2)}$\\
$^{(1)}$ J.A.D. Laboratory, University of Nice Sophia-Antipolis, France. hamad.hazim@unice.fr, br@unice.fr\\
$^{(2)}$ ISVR,  Southampton University, UK. nsf@isvr.soton.ac.uk\\
In collaboration with  Thales Alenia Space, France}

\maketitle
\section{ABSTRACT}
The mass reduction of satellite solar arrays results in significant panel flexibility, so possibly striking one another dynamically leading ultimately to structural damage. To prevent this, rubber snubbers are
mounted at well chosen points of the structure and they act as one sided linear
spring; as a negative consequence, the dynamic of these panels becomes nonlinear.
 The finite element approximation is used to solve partial differential equations governing the structural
dynamic. The models are validated and adjusted with experiments done in the ISVR laboratory, Southampton university.
\begin{figure}[hbtp]
\begin{center}
\includegraphics[width=3cm]{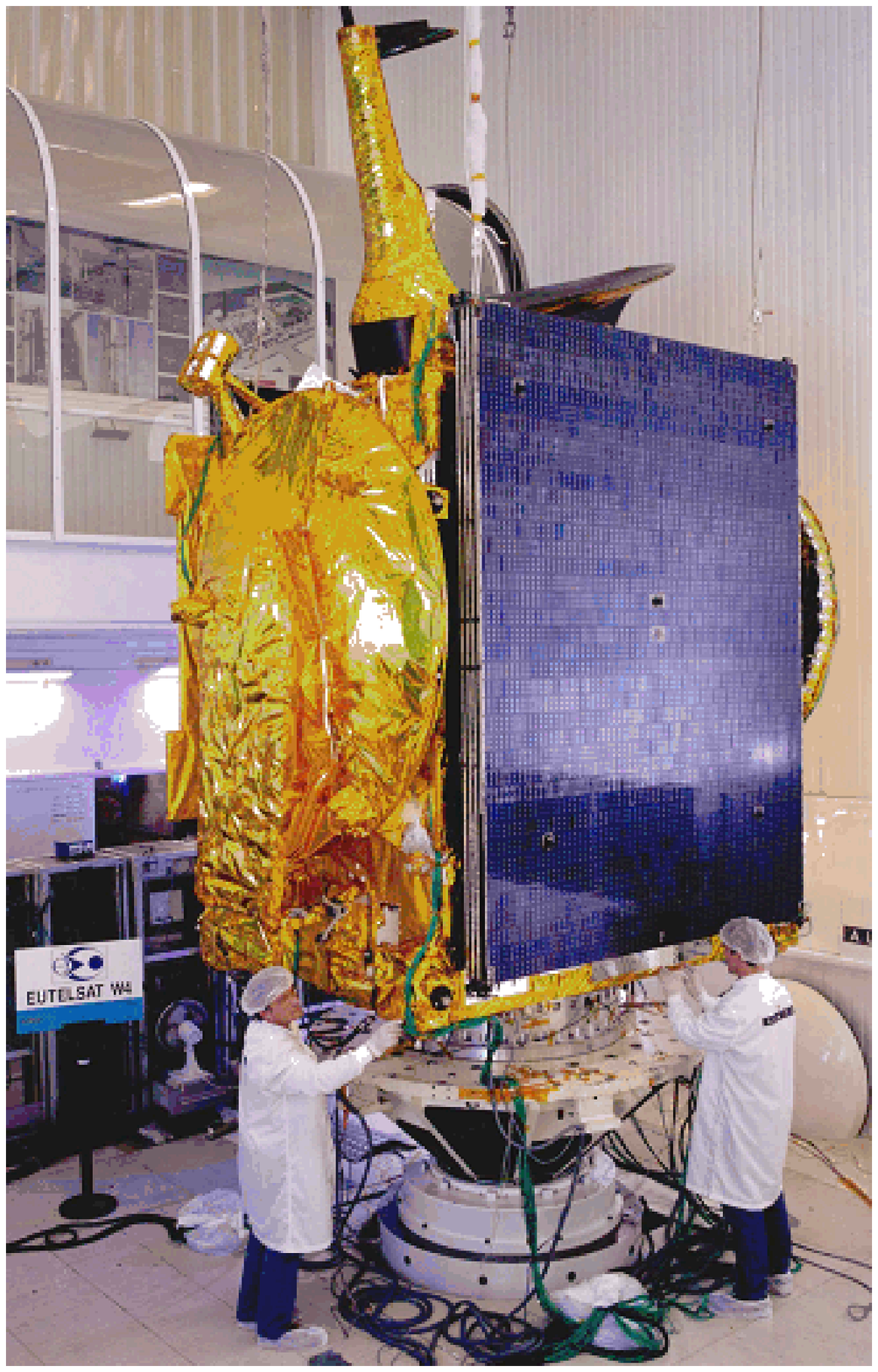}
\hspace{0.1cm}
\includegraphics[width=3.5cm]{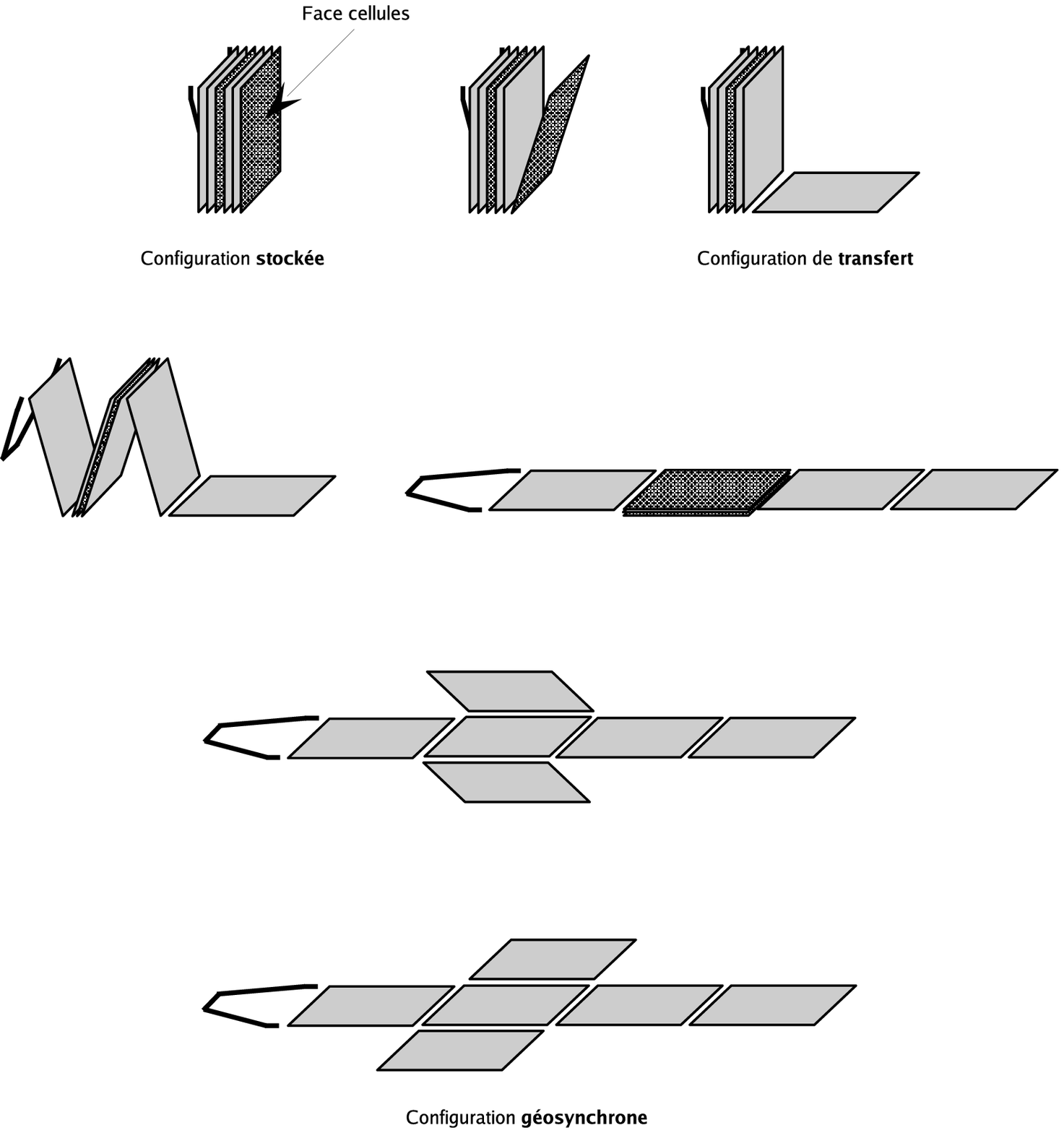}
\end{center}
\caption{\small{Solar array of the satellite under a test on a shaker}} \label{1}
\end{figure}
\section{INTRODUCTION}
\begin{figure}[hbtp]
\begin{center}
\includegraphics[width=7cm,height=4cm]{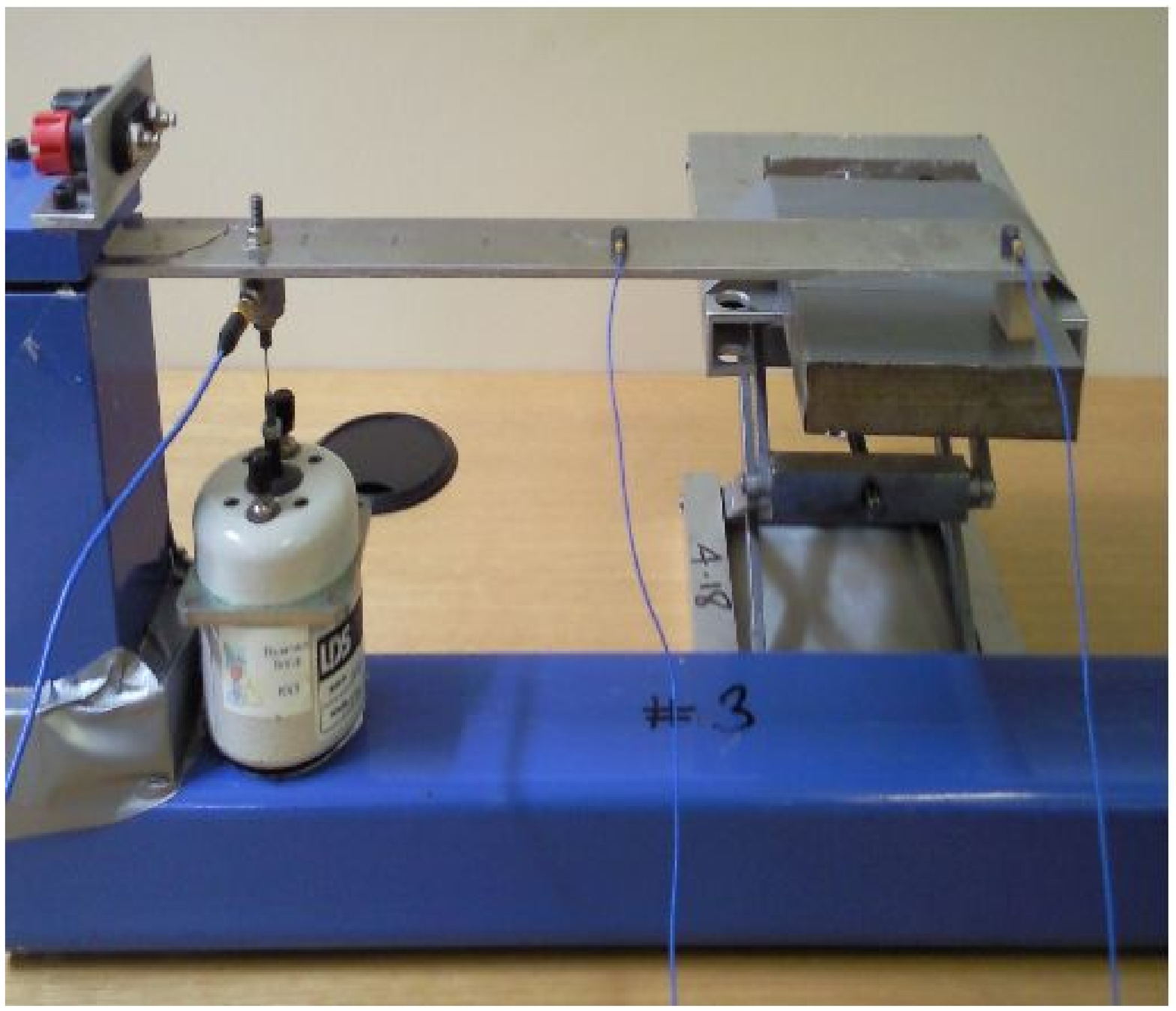}

\includegraphics[width=7cm,height=3.25cm]{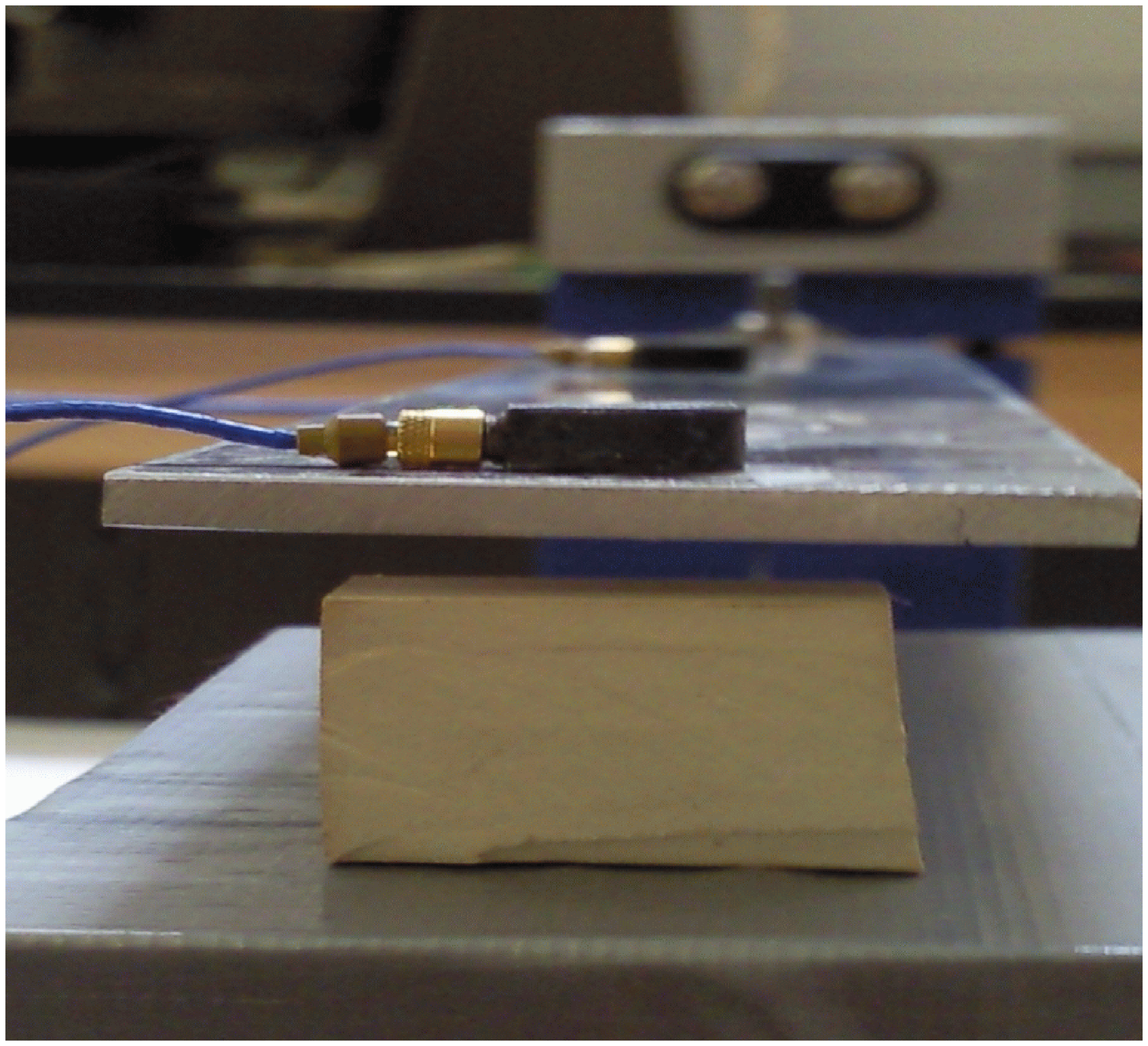}\label{rig1}
\end{center}
\caption{\small{The rig used for the experiments: a linear clamped-free beam in contact with a rubber}} \label{1}
\end{figure}

The study of the total  dynamic behavior of solar arrays in a folded position with snubbers  are so complicated, that to simplify, a solar array is modeled by a clamped-free Bernoulli beam with one-sided linear spring. This system is subjected to a periodic excitation force.\\
The real state of the problem is close to a beam with unilateral contact subjected to a periodic imposed displacement of the base, but the
dynamical behavior of the system do not change  significantly if the imposed displacement is replaced by a periodic force excitation, the configuration used
is easiest to be realized from a technical point of view as the rig is very simple to built [see figure $\ref{rig1}$]\\
The present study is to simulate the behavior of a beam which strikes a snubber under a periodic excitation. As the interest is to deal with the first three eigen frequencies, the beam is modeled by ten finite elements.\\
The results are presented and compared in the frequency domain, the FFT is applied to the predicted and the experimental displacement of the last node of the beam finite element. The mass effect of the force transducer is taken in account in the finite element code. Note that no treatment is done or used from the acquisition system as the problem is nonlinear, the transfer function and the other functions use linear assumptions. The time signal is taken and the processing is done using an external software (Scilab $\cite{a6}$). The numerical predictions compared to experimental results show very good agreement.
\newpage
\section{THE MODEL}
\begin{figure}[hbtp]
\begin{center}
\includegraphics[width=8cm,height=4cm]{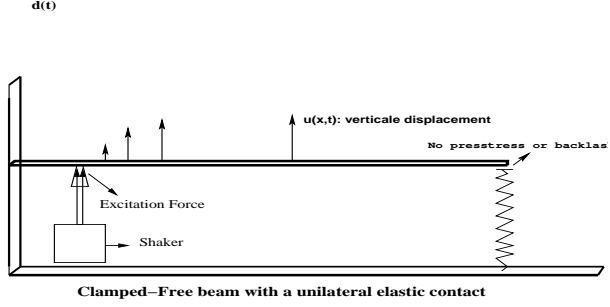}
\end{center}
\caption{\small{beam system with an unilateral spring under a periodic excitation}} \label{1}
\end{figure}
The beam motion with a snubber can be modeled as:
\begin{equation}
\rho S\ddot{u}(x,t)+EIu^{(iv)}(x,t)=F(t)\delta_{x_0}-k_r u(x_1,t)_-
\end{equation}
$u(0,t)=0$, $\partial_x u(0,t)=0$
\begin{equation}
u(x,t)_-=\left\lbrace
\begin{array}{rl}
u(x,t) & \ if \ u \leq 0\\
0  & \ if \ u > 0\\
\end{array} \right.
\end{equation}
The classical Hermite cubic finite element approximation is used, it yields an ordinary differential system in the form:
\begin{equation}
M\ddot{q}+Kq=-[k_{r}(q_{n_1})_-]e_{n_1}+F(t) e_{n_2}
\end{equation}
Where $M$ and $K$ are respectively the mass and stiffness assembled matrices, q is the vector of degrees of freedom of the beam,
$q_i=(u_i,\partial_xu_i)$, $i=1,...,n$, where $n$ is the size of $M$, $n_1$ and $n_2$ are the indices of the nodes where the spring and the excitation force meet the beam respectively. Numerical time integration was performed using ODE numerical integration for 'stiff' problems, package ODEPACK is called and it uses the BDF method.\\
The use of a small electrodynamic shaker yields a technical problem due to the reaction of the beam, the input force $F(t)$ can not be a simple sine wave. To deal with this problem, a force transducer is fixed between the shaker and the beam. The numerical codes use the actual force time signal coming from the acquisition system which is periodic. Figure $[\ref{f1}]$ shows the input force signals with its spectrum contents.\\

\section{EXPERIMENTAL AND NUMERICAL RESULTS}
\begin{figure}[hbtp]
\begin{center}
\includegraphics[width=8cm,height=4.05cm]{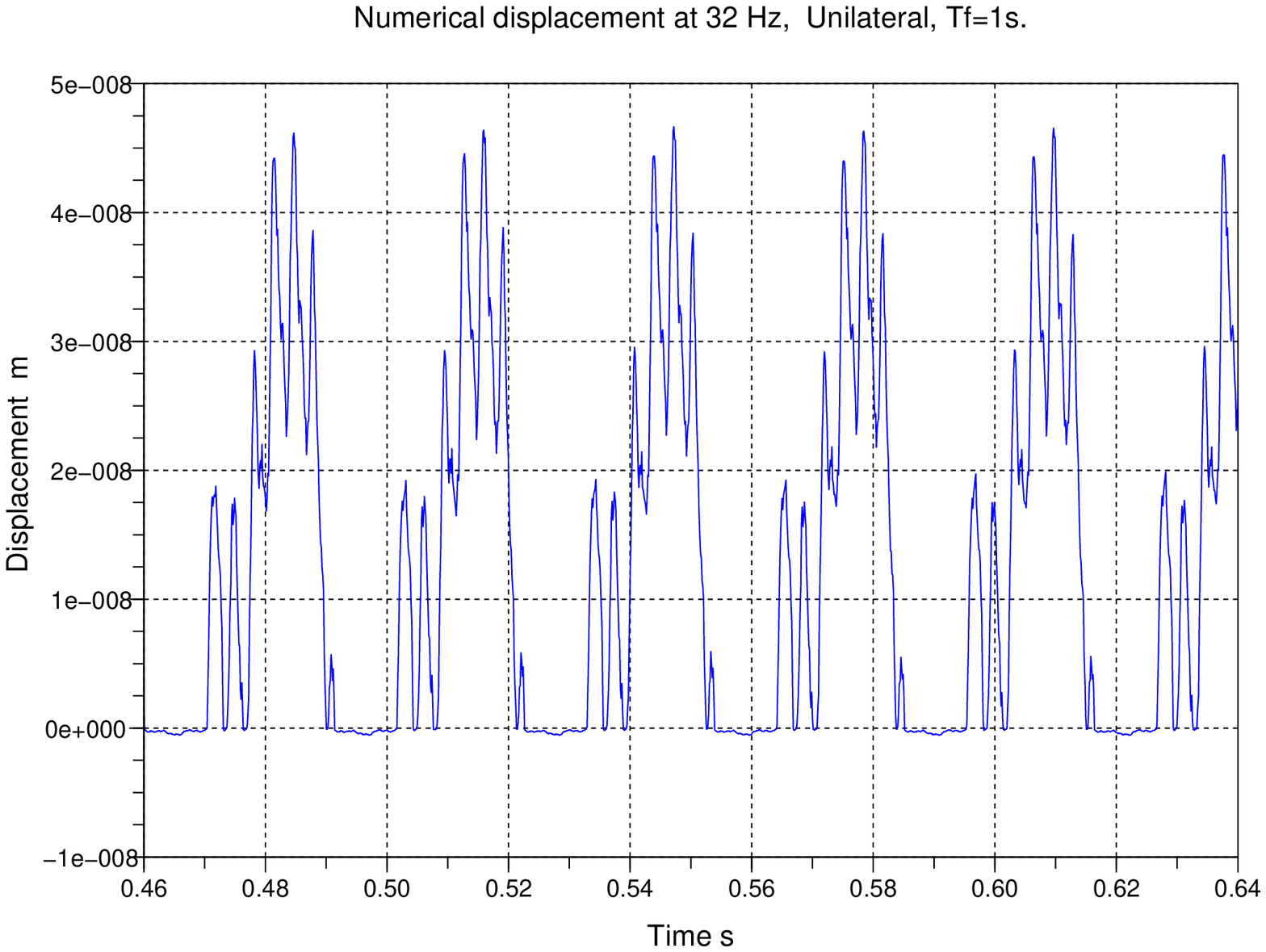}
\includegraphics[width=8cm,height=4.05cm]{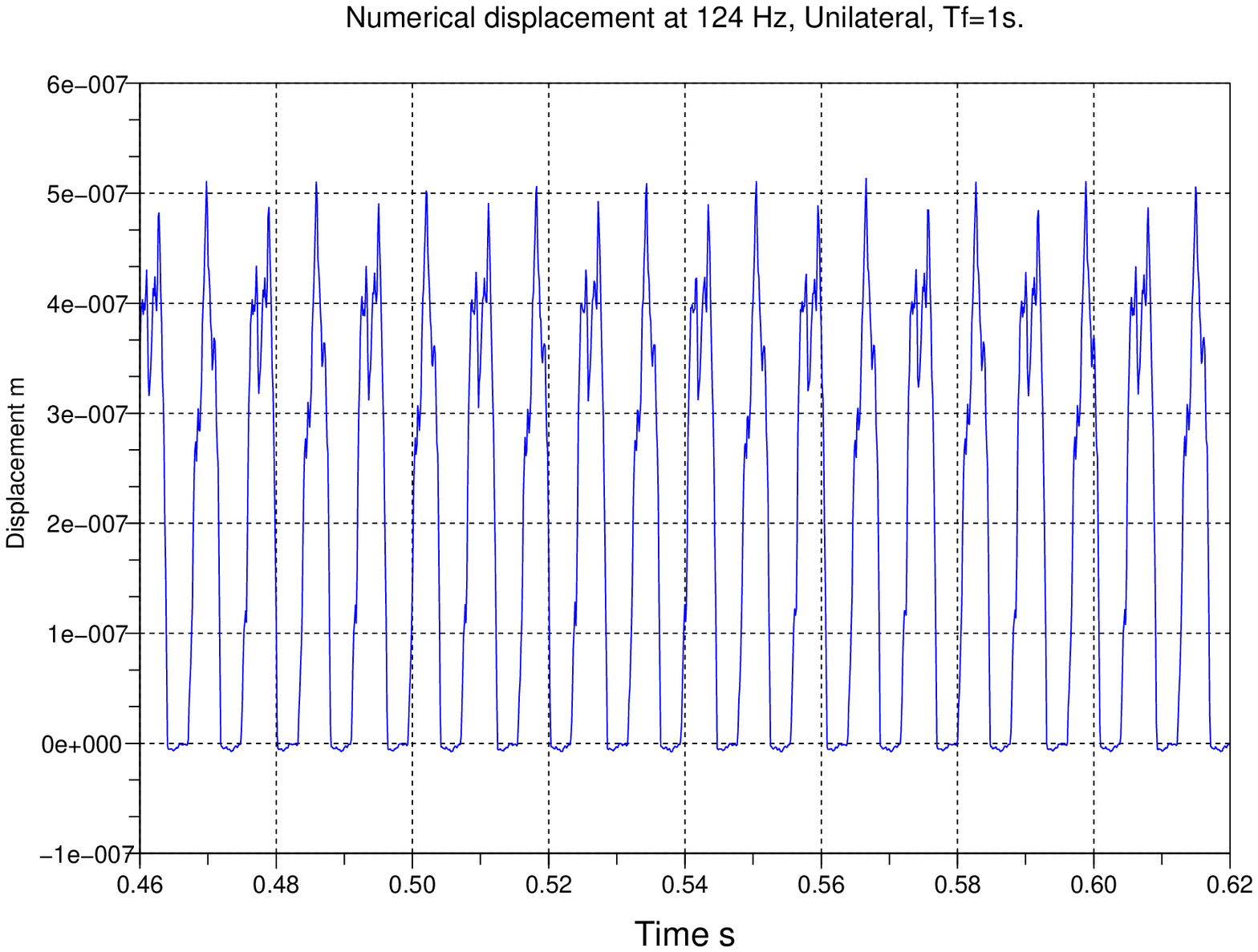}
\caption{The predicted displacements for an excitation of $32 \;Hz$ and $124\; Hz$, the high stiffness of the rubber yields almost  a positive displacement}\label{d}
\end{center}
\end{figure}

\begin{figure}[hbtp]
\includegraphics[width=8cm,height=11.5cm]{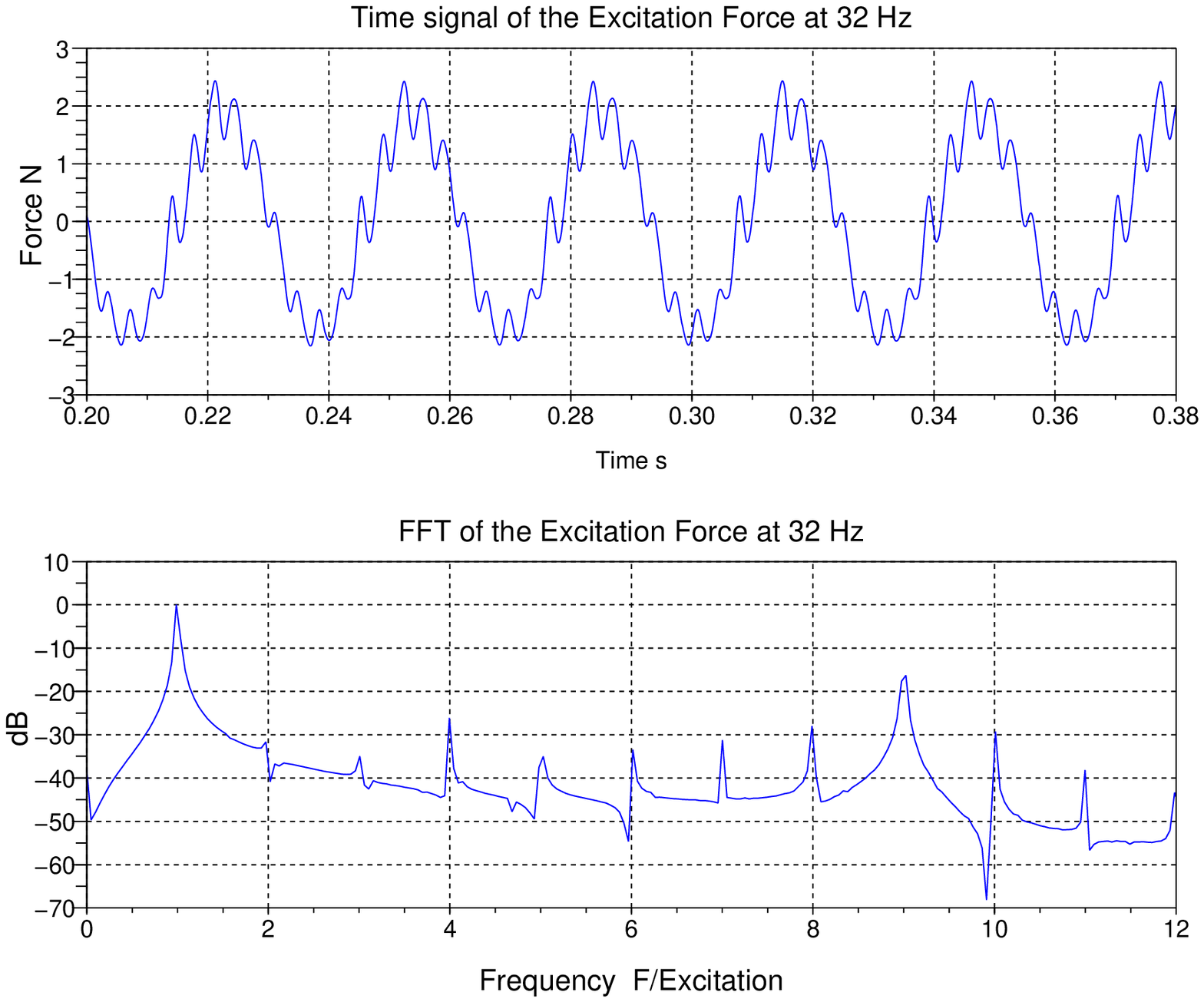}
\includegraphics[width=8cm,height=11.5cm]{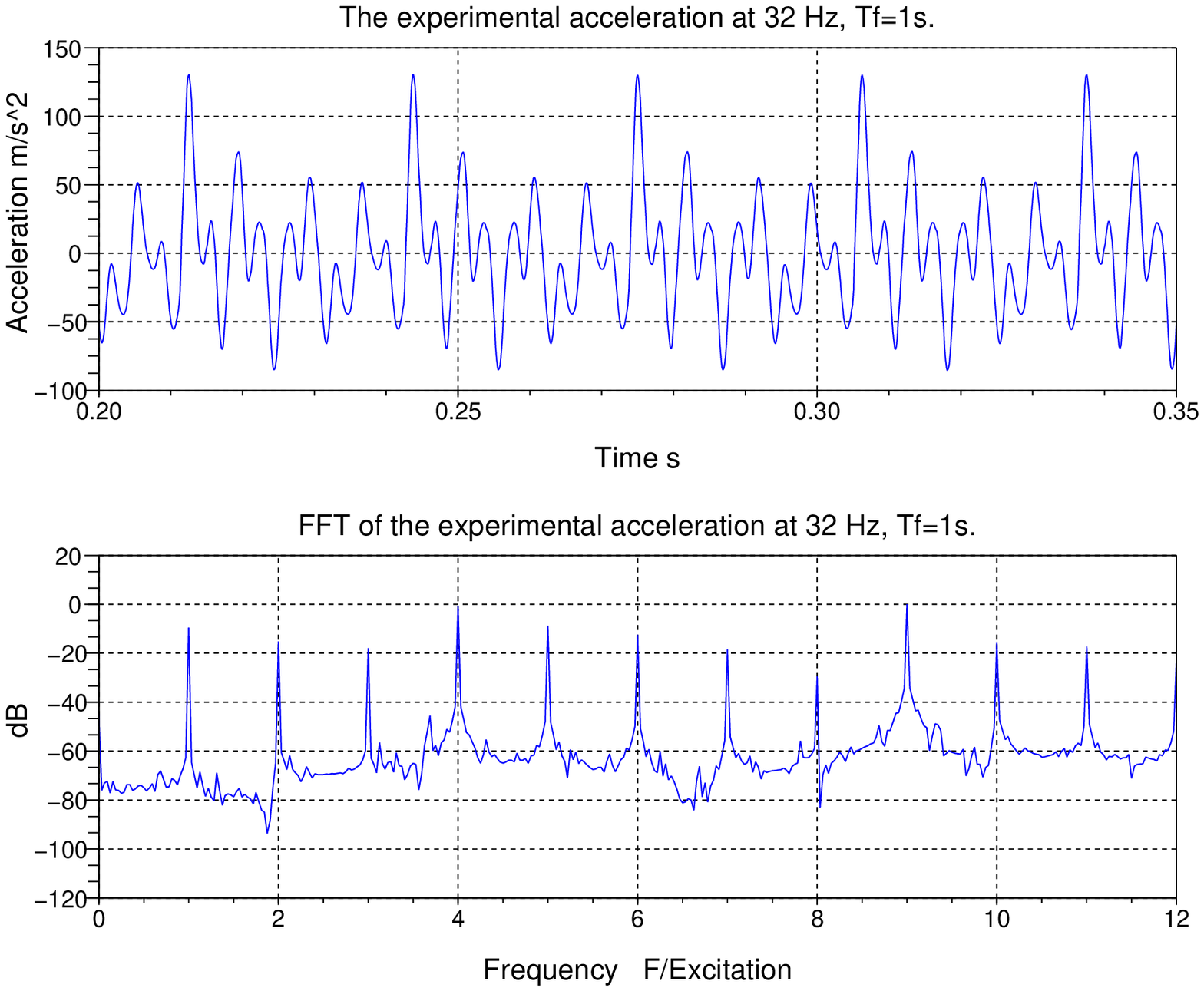}
\caption{The time signal and its FFT of the input force and the experimental acceleration for an excitation of $32 Hz$}\label{f1}
\end{figure}

\begin{figure}[hbtp]
\includegraphics[width=8cm,height=11.5cm]{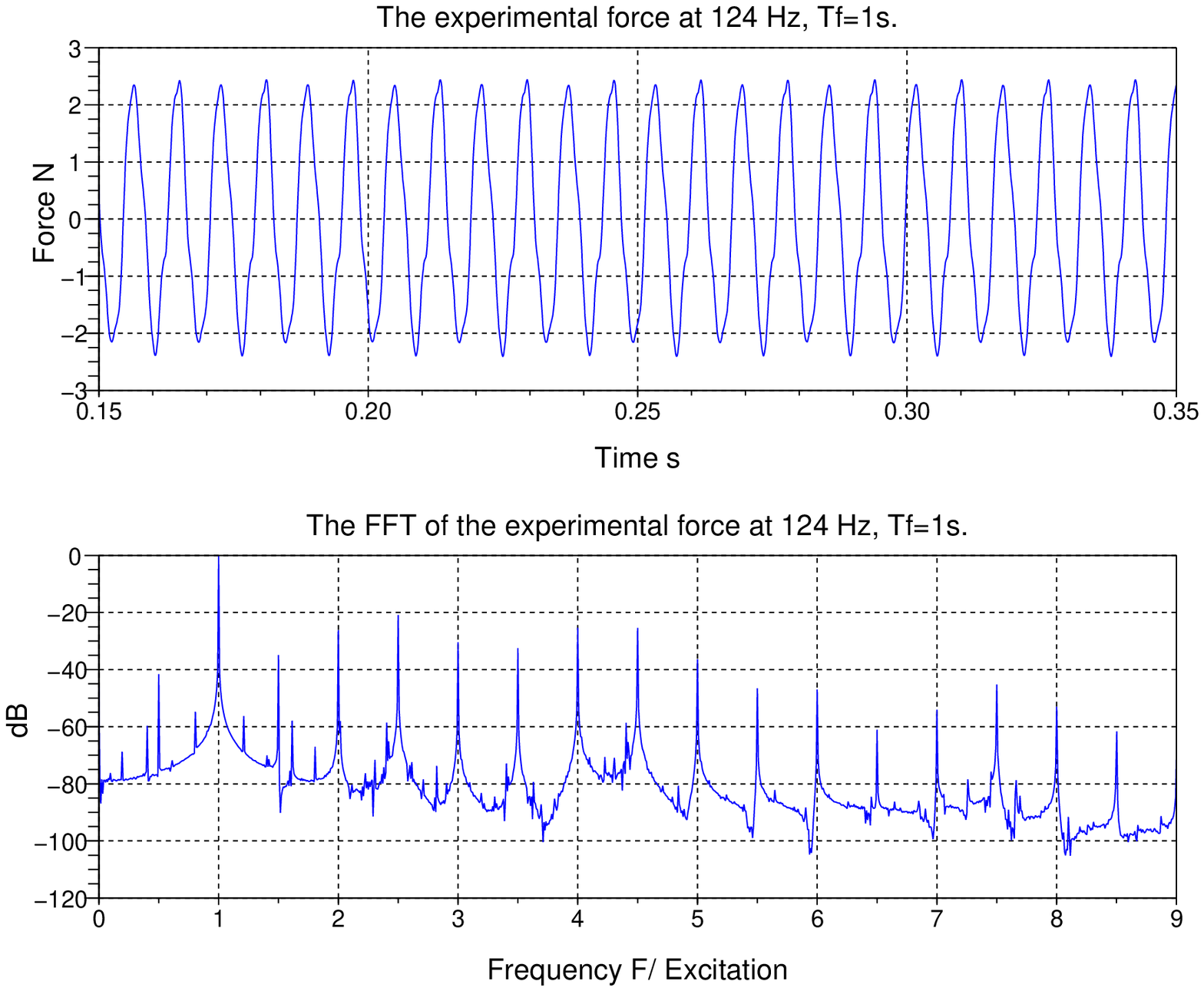}
\includegraphics[width=8cm,height=11.5cm]{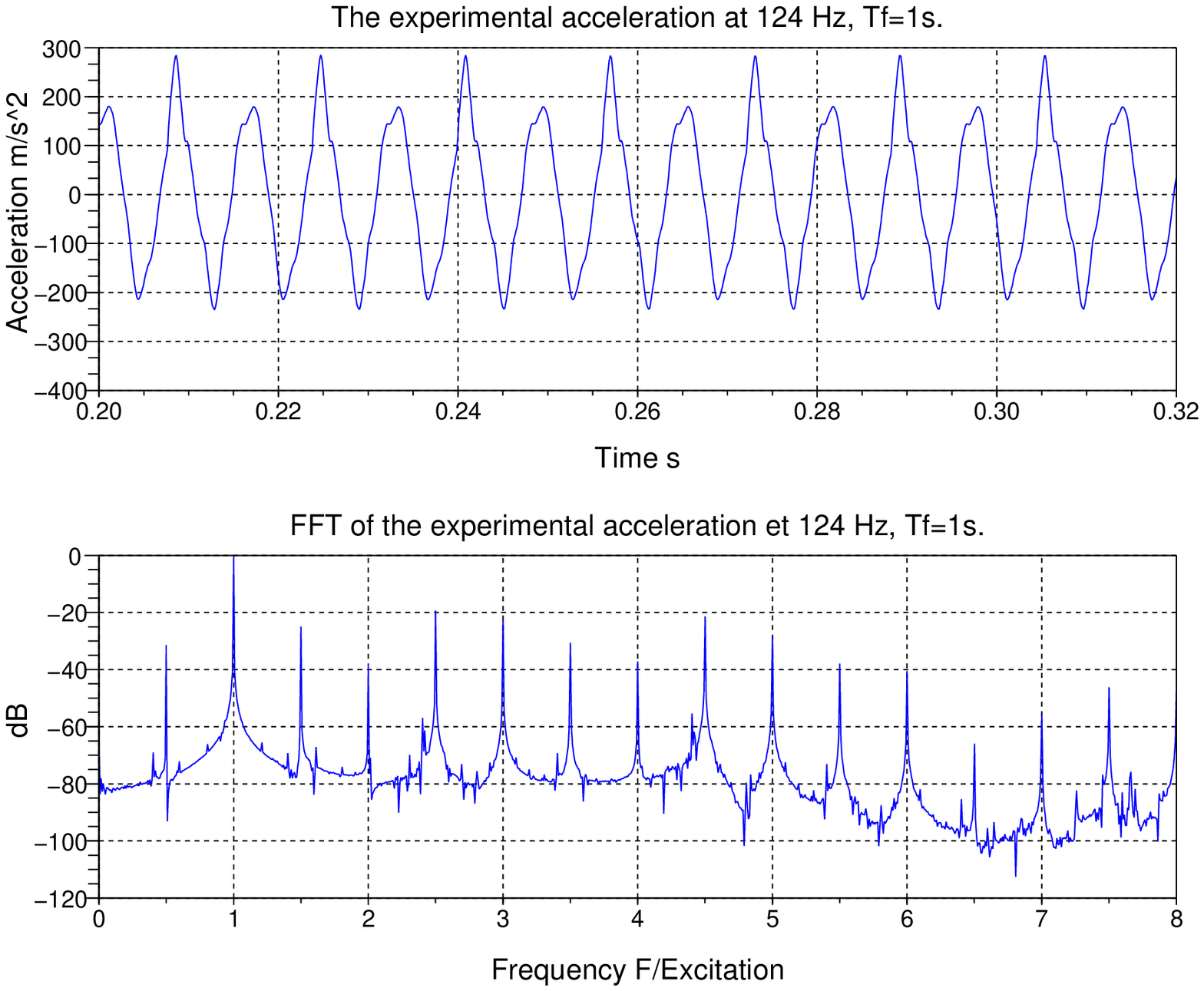}
\caption{The time signal and its FFT of the input force and the experimental acceleration for an excitation of $124Hz$}\label{f2}
\end{figure}

\begin{figure}[hbtp]
\includegraphics[width=9cm,height=8cm]{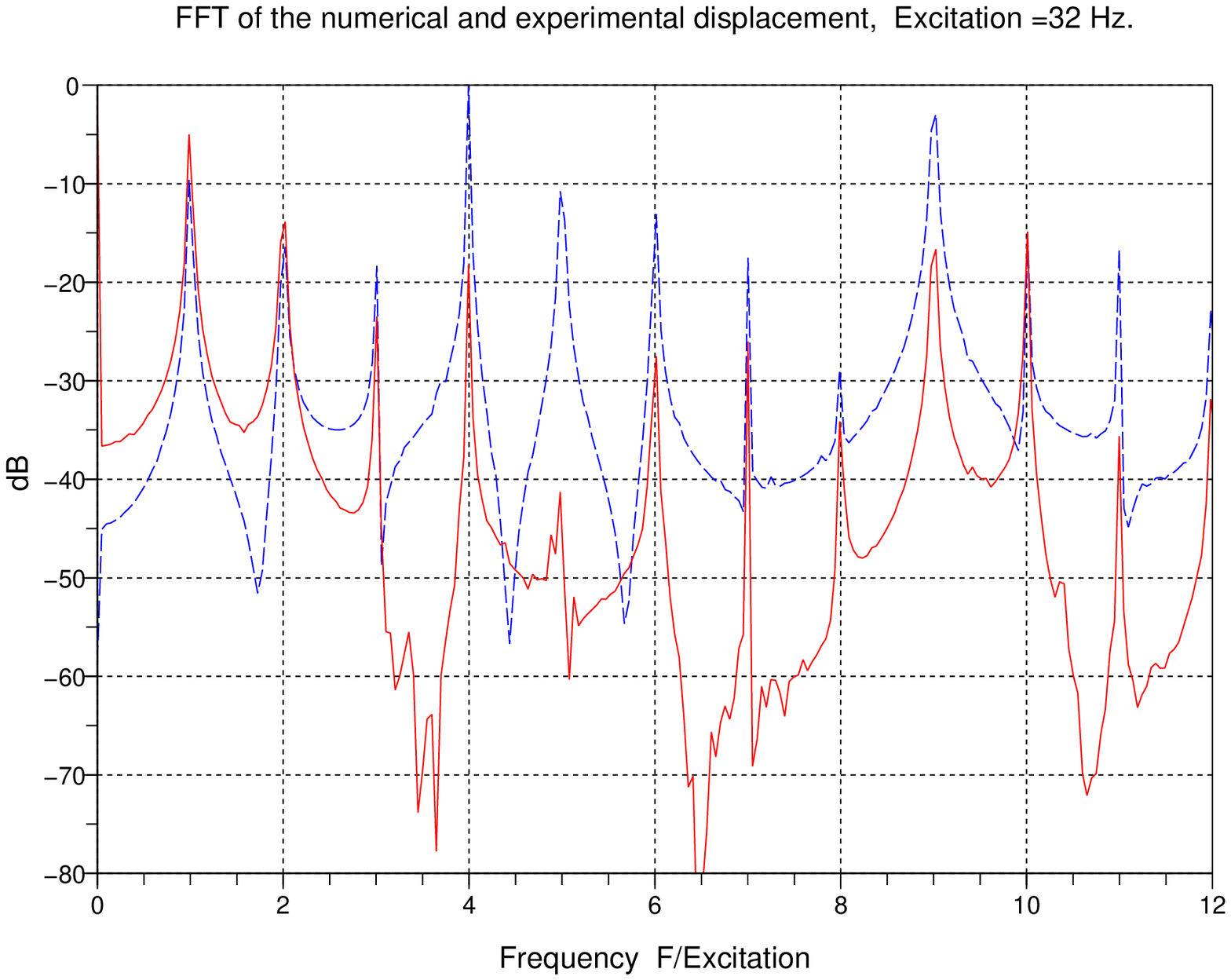}
\includegraphics[width=9cm,height=8cm]{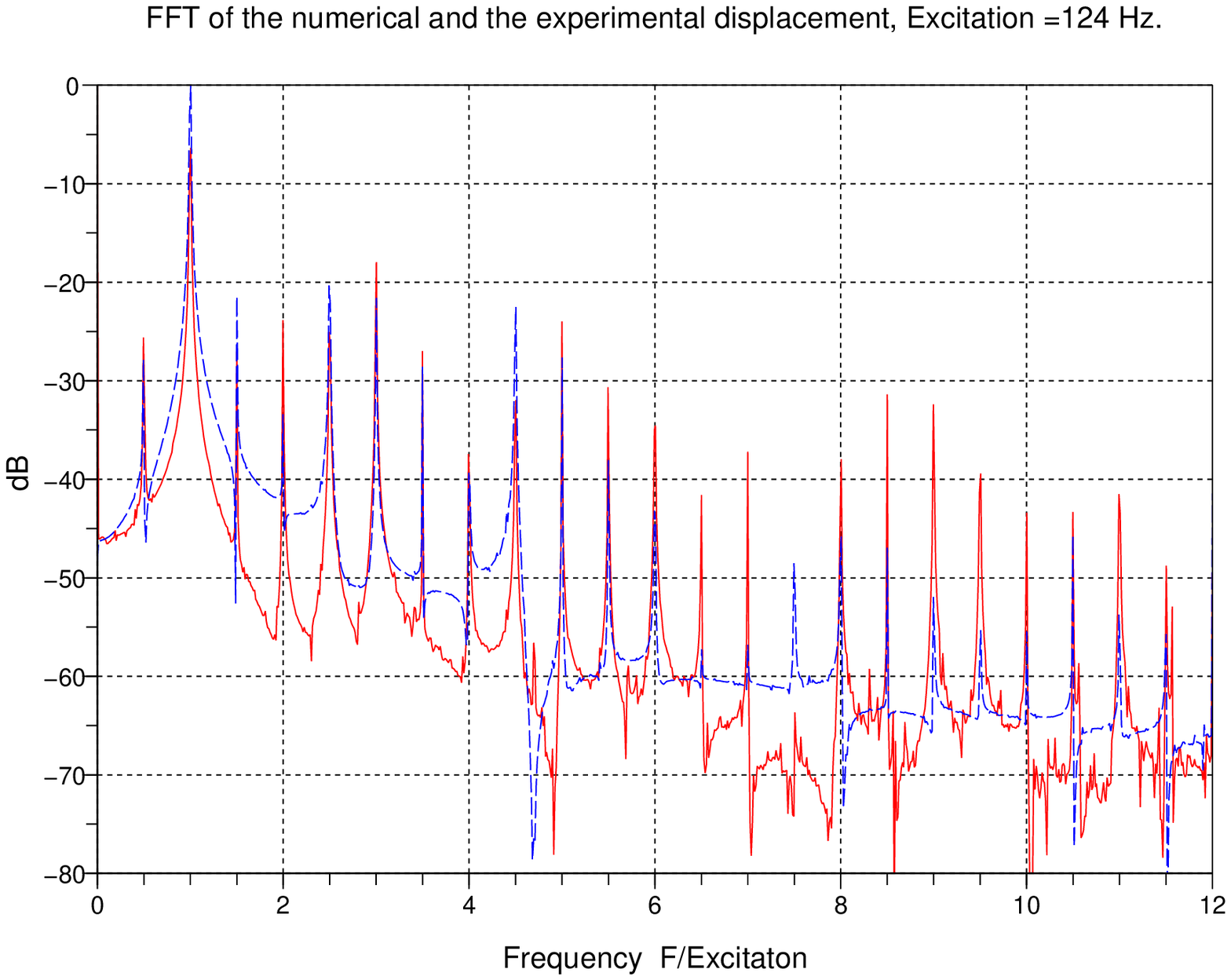}
\caption{dashed line: FFT  and the experimental displacements, solid line: FFT of the numerical displacements. Two different excitation: $32 Hz$ and $124 Hz$}\label{fft}
\end{figure}

Figures [$\ref{rig1}$] and [$\ref{1}$] give a clear idea of the experiment setup. In the table, all the parameters used in the simulations and the experiments are shown as well as the  first three predicted eigen frequencies of the linear beam without the spring. The integration time is fixed at $t=1s$ for all the simulations and the experiment sequences which is five times the greater period of the system.\\
Figure $[\ref{fft}]$ shows the FFT of the numerical and the experimental displacements for $32 HZ$ and $124 Hz$, the height of the peaks are normalized by the maximum; the predicted frequencies found are exactly the same measured for a large number of harmonics. However, a small shift in the height of these pikes appears from the fifth harmonic, this is maybe due to the low number of finite elements used to model the beam; some investigations are in progress to understand this aspect. Other tests with random excitations have shown good agreement. The effect of the unilateral contact is underlined, the input frequency is split into its all harmonics.\\
Figures  $[\ref{f1}]$ and  $[\ref{f2}]$ show the input excitation force, it is clear from the time signal and from its frequency content that this force is not a single sine wave, but it is periodic.\\
Figure $[\ref{d}]$ shows the predicted displacement for an excitation of $32 Hz$ and $124 Hz$, the displacement is almost positive, this is due to the high stiffness of the spring but the time response is still periodic.\\
Note that the frequencies $32 Hz$ and $124Hz$ are chosen and presented herein because they represent the first two nonlinear normal modes of the nonlinear normal mode, they are computed numerically by a numerical sweep test. Non linear normal modes $NNM$ of a nonlinear differential systems are an extension of the well known linear normal modes,another interest of the authors with some papers in progress.
\tiny{
\begin{equation}\label{tab1}
\begin{tabular}{|c|c|c|c|c|c|c|}
\hline
& & & \\
\textbf{beam length} & \textbf{beam width} & \textbf{beam thickness}& \textbf{beam's Young modulus} \\
& & & \\
\hline
& & & \\
$0.35 m $& $0.0385 m$& $0.003m$ & $69\times10^9 N/m^2$\\
& & & \\
\hline
&  &  & \\
\textbf{Spring stiffness} & $1^{"}$ \textbf{linear eigen frequency} & $2^{"}$ \textbf{linear eigen frequency} & $3^{"}$\textbf{ linear eigen frequency} \\
 &  &  & \\
\hline
& & & \\
 57.14 $KN/m$&$19.97 Hz$ & $122.2 Hz$ & $318.8Hz$ \\
&  & & \\
\hline
\end{tabular}
\end{equation}
}

\vspace{7cm}



\begin{thebibliography}{1}
\small{
\bibitem[1]{r1} J. H. Bonsel, R. H. B. Fey. and H. Nijmeirjer: Application of a Dynamic Vibration Absorber to a Piecewise Linear Beam System.
\textit{  Nonlinear Dynamics} 37: 227-243, 2004.

\bibitem[3]{a4} D.Jiang,C.Pierre, S.W. Shaw: Large-amplitude non-lin\'{e}ar normal modes of
piecewise linear systems. \textit{Journal of sound and vibration}
272 (2004) 869-891

\bibitem[3]{r2} Rob H.B. Fey and  Berend Winter, Jaap J. Wijker\\Sine sweep and steady - state response of a simplified solar array model with nonlinear support.
DETC99/VIB-8027, 1999 ASME Design Engineering Technical Conferences,
Las Vegas, Nevada, USA.

\bibitem[4]{r4} R.H.B. Fey and F.P.H. van Liempt :\\ Sine sweep and steady-state response of a simplified solar array model with nonlinear elements.
International conference on Structural Dynamics Modelling, 3-5 june
2002, pp.201-210.

\bibitem[5]{G} M G\'{e}radin Th\'{e}orie des vibrations: applications \`{a} la dynamique des structures. Masson Paris 2eme edition 1996

\bibitem[6]{JB}  St\'{e}phane Junca et Bernard Rousselet: Asymptotic Expansions of Vibrations with Unilateral Contact, proceeding of the GDR 2501-Anglet 2-6 Juin 2008.

\bibitem[7]{r3} Richard H. Rand:Lecture Notes on Nonlinear Vibrations, Dept. Theoretical Applied Mechanics Cornell University Ithaca NY 14853
\bibitem[8]{a6} Scilab [www.scilab.org]
\bibitem[9]{r5} E.L.B. Van De Vorst, M.F.Heertjes\ D.H. Van Campen and A. De Kraker
and R.H. B. Fey: Experimental and numerical analysis of the steady
state begavior of a beam system with impact. \textit{Journal of
Sound and Vibration} 1998 212(2), 210-225}
\end{thebibliography}
\end{document}